\documentclass[twocolumn, 10pt]{article}

\usepackage{algorithmic}
\usepackage{url}
\usepackage{tikz}
\usepackage{float}
\usepackage{graphicx}
\usepackage{blkarray}
\usetikzlibrary{arrows.meta}
\usepackage{caption}
\usepackage{amsmath}
\usepackage{amssymb}
\usepackage{multirow}
\usepackage{booktabs}
\usepackage{cite}
\usepackage[margin=0.5in]{geometry}

\newtheorem{proposition}{Proposition}[section]

\DeclareMathOperator{\hd}{d_{H}}
\DeclareMathOperator{\sv}{sv}

\begin{document}
\title{A Supervised Learning Approach to Rankability}
\author{Nathan McJames$^{1,2,*}$, David Malone$^{1,2}$, and Oliver Mason$^{2}$}
\date{%
    $^1$Hamilton Institute, Maynooth University, Maynooth, Co. Kildare, Ireland\\%
    $^2$Department of Mathematics and Statistics, Maynooth University, Maynooth, Co. Kildare, Ireland\\%
    $^*$Corresponding author. E-mail address: nathan.mcjames.2016@mumail.ie\\[2ex]%
    \today
}

\twocolumn[
  \begin{@twocolumnfalse}
    \maketitle
    \begin{abstract}
    \noindent The rankability of data is a recently proposed problem that considers the ability of a dataset, represented as a graph, to produce a meaningful ranking of the items it contains. To study this concept, a number of rankability measures have recently been proposed, based on comparisons to a complete dominance graph via combinatorial and linear algebraic methods. In this paper, we review these measures and highlight some questions to which they give rise before going on to propose new methods to assess rankability, which are amenable to efficient estimation. Finally, we compare these measures by applying them to both synthetic and real-life sports data.
    \end{abstract}
  \end{@twocolumnfalse}
  ]

\section{Introduction}
The rankability of data is a new problem first posed by the authors of \cite{edgeR}, which they describe as ``a dataset's inherent ability to produce a meaningful ranking of its items''.  Questions of ranking and rankability arise naturally in connection with datasets represented as directed graphs.  The graphs may be weighted, with the graph structure representing the pairwise comparisons between the elements of the dataset, which are the nodes/vertices of the graph.  Before attempting to obtain a ranking of the elements in a dataset, it first makes sense to consider how suitable the dataset is for this task.  Intuitively, certain structures seem more ``rankable'' than others. For example, we would expect trees to be more rankable than cycles or complete graphs.  The ideal situation, from the point of view of rankability, is where there is a clear unambiguous ordering of the nodes of the graph.  This is described by a \emph{complete dominance graph}, where there is some ordering $i_1, i_2, \ldots, i_n$ of the nodes such that the edges of the graph are given by $\{(i_p, i_q): p < q\}$.  

To quantify the concept of rankability, the authors of \cite{edgeR} defined a new measure based on the number of edge additions or removals that are necessary to convert a directed graph into a complete dominance graph.  This first measure applies to unweighted graphs only.  More recently, the authors of \cite{specR} proposed a second rankability measure based on the spectral variation of the Laplacian matrix of a complete dominance graph. This measure is much more widely applicable, and can be used for weighted graphs, but struggles to return meaningful results for sparse data. In both cases, the authors focus on the \emph{rankability} of the data rather than \emph{ranking} it.

In this paper we develop a framework for rankability based on the generation of graphs with a chosen/candidate \emph{target rankability}. We give an example of this by making perturbations to complete dominance graphs. We will show that this approach enables evaluation of the performance of previous rankability measures. It also enables the use of a supervised learning approach for the measurement of rankability. We will demonstrate this approach by training a random forest on synthetic data generated with a chosen target rankability. We will also test random forest regression models on graphs generated with a different graph generation model and real life sports data. Results from this experiment suggest that the random forest is capable of performing well on the out of sample synthetic data and the real-life sports data also.

\subsection{Ranking}
The ranking of data is a key data science task which has many important applications. Search engines rank web pages by their relevance to a given search query \cite{searchEngines}, and media companies rank films and television series by their similarity to a viewer’s likes and dislikes \cite{netflix}. A related application is the ``Analytic Hierarchy Process'' in Decision Science \cite{AHP}, where pairwise comparisons between alternatives generate \emph{symmetrically reciprocal matrices} and these are then used to produce a transitive or linear ordering.  

The ranking of data is commonly formulated in terms of the Linear Ordering Problem (LOP), which has been the subject of much research in combinatorial optimization and related applications \cite{LOP1, LOP2}.  The objective of the LOP is as follows: Given information on some or all of the pairwise comparisons of a set of $n$ items, produce the linear ordering of these items which is in the strongest agreement with the available data. This optimal linear ordering of the $n$ items is found by solving an integer programming problem \cite{LOP3}.

One disadvantage of ordering/ranking items with the LOP is that the associated integer programming problem may not have a unique solution. A second disadvantage is that the ordering returned by the algorithm may not be \emph{unquestionable}. That is, there may be other orderings with only marginally worse values of the objective function. In the simplest case, the best ordering of the items in a dataset will be immediately obvious. This occurs when the relationships within a dataset form a complete dominance graph. A complete dominance graph provides us with a \emph{unique and unquestionable ranking}. This question of how unique and how unquestionable a given ordering of the items within a dataset is is closely related to rankability.

\subsection{Rankability}
The problem of rankability of data was only recently posed by the authors of \cite{edgeR}. This topic can be described as an assessment of the appropriateness of creating a linear ordering of the items in a dataset. This is an important assessment to make because a ranking algorithm applied to a dataset which is not rankable may yield misleading results. In their paper they define a rankability measure which we will refer to as $R_{e}(G)$, and outline three key requirements of a good rankability measure:

\begin{enumerate}
  \item It should be effective. In other words it should return a result that makes sense. It should assign a high rankability to a complete dominance graph and low rankability to a graph with no way of creating a sensible ordering such as a completely connected graph.
  \item It should be efficient, i.e. be able to compute the rankability of a graph in a reasonable amount of time.
  \item It should be algorithm agnostic. The rankability measure should be a function solely of the graph and nothing else. It is not a measure of the quality of any ordering of the data.
\end{enumerate}

Requirement one above is particularly relevant to the content of this paper. In the case of small graphs, one can examine the edges by eye, and decide approximately how rankable the graph is. In this way it is also possible to evaluate the efficacy of a rankability measure by checking that it returns sensible results for small graphs.  For larger graphs, however, this is more difficult, and a rankability measure that appears to perform well on small graphs might yield misleading results on larger graphs or vice versa. Later we will investigate this problem in more detail and suggest a method for evaluating a rankability measure's performance based on simulations involving synthetic data.

More recently, the authors of \cite{specR} suggested a second rankability measure, which we will refer to as $R_{s}(G)$, and which is designed to address requirement two above which is not fulfilled by $R_{e}(G)$. This rankability measure is much more widely applicable, but struggles to return meaningful results where not all pairwise comparisons have been considered, so fails to fulfil requirement number one in such cases. In their conclusion, they suggest future research should attempt to address this issue, and also seek to define a rankability measure based on different graph properties. 

In this paper, we attempt to provide a solution to all three of the above requirements by tackling the problem of rankability with a supervised learning approach. We begin by recalling the previously defined rankability measures, $R_{e}(G)$ and $R_{s}(G)$, and highlight some of their properties. We then suggest a method for the evaluation of a rankability measure's performance which is based on generating graphs with a known target rankability. We also show that by generating graphs with a known target rankability, it is possible to train a random forest regression model to measure the rankability of graphs. We support this supervised learning approach by applying it to out of sample synthetic data and verifying that it still yields sensible results. We also apply the rankability measures to real world sports data and observe a strong correlation between them, suggesting a level of mutual agreement.

\section{Observations on Previous Rankability Measures}
\label{sec:previous}

\subsection{Edge Rankability}
\label{EdgeR}
Proposed by the authors of \cite{edgeR}, $R_{e}(G)$ is a rankability measure based on the smallest number of edge additions or removals that are necessary to convert a directed graph into a complete dominance graph. This smallest such number is denoted by $k$. It also takes into account the number of different complete dominance graphs that are reachable by making just $k$ changes. This second number is denoted by $p$. Having determined these two values, the graph's rankability is defined as 
\begin{equation}\label{eq:Re}
R_{e}(G)=1-\frac{kp}{k_{max}p_{max}},
\end{equation}
where $k_{max}$ and $p_{max}$ are normalisation factors to keep the resulting rankability bounded between 0 and 1. They represent the largest possible values that $k$ and $p$ can take, which in the case of an $n$ vertex graph are $\frac{n(n-1)}{2}$ and $n!$ respectively.

Recall that the adjacency matrix of an unweighted directed graph $G$ with vertices given by $\{1, \ldots, n\}$ is the matrix $A$ in $\mathbb{R}^{n \times n}$ with $a_{ij} = 1$ if $(i,j)$ is an edge of $G$ and $a_{ij} = 0$ otherwise. Accordingly, the adjacency matrices $D1$ and $D2$ below correspond to two different orderings of the vertices of the graph in Figure 1.

\[D1=
\begin{blockarray}{ccccc}
&1&3&4&2\\
\begin{block}{c(cccc)}
  1&0 & 1 & 1 & 1 &\\
  3&0 & 0 & 0 & 1 &\\
  4&0 & 0 & 0& 1 &\\
  2&0 & 0 & 0 & 0 &\\
\end{block}
\end{blockarray}
\ \ \ D2=
\begin{blockarray}{ccccc}
&1&4&3&2\\
\begin{block}{c(cccc)}
  1&0 & 1 & 1 &1&\\
  4&0 & 0 & 0 &1&\\
  3&0 & 0 & 0&1&\\
  2&0 & 0 & 0 &0&\\
\end{block}
\end{blockarray}
\]

The graph in Figure 1 is nearly a complete dominance graph.
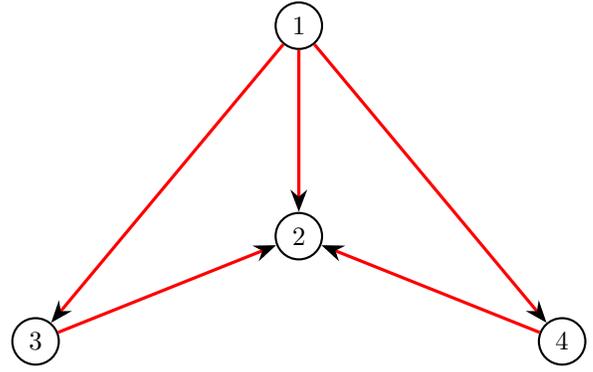
\begin{figure}
\begin{center}
\begin{tikzpicture}[scale=0.7]
\begin{scope}[every node/.style={circle,thick,draw}]
    \node (1) at (0,0) {1};
    \node (2) at (0,-4) {2};
    \node (3) at (-5,-6) {3};
    \node (4) at (5,-6) {4};
\end{scope}
\begin{scope}[>={Stealth[black]},
              every node/.style={fill=white,circle},
              every edge/.style={draw=red,very thick}]
    \path [->] (4) edge (2);
    \path [->] (3) edge (2);
    \path [->] (1) edge (2);
 \path [->] (1) edge (3);
 \path [->] (1) edge (4);
\end{scope}
\end{tikzpicture}
\caption{Example graph for calculating $R_{e}(G)$ and $R_{s}(G)$}
\label{Example graph for calculating edgeR}
\end{center}
\end{figure}
Vertex 1 is clearly the best and vertex 2 is clearly the worst. The only indecision lies in the comparison between vertices 3 and 4. This uncertainty is reflected in adjacency matrices $D1$ and $D2$, both of which require an edge addition to achieve the characteristic upper triangular pattern of adjacency matrices for complete dominance graphs. Matrix $D1$ requires an edge addition at (3, 4), and $D2$ requires an edge addition at (4, 3). Both adjacency matrices require just one edge change, and therefore $k$ for this graph must be equal to 1. Since $D1$ and $D2$ are the only adjacency matrices requiring just one change, $p$ must be equal to 2. For a graph with 4 vertices, $k_{max}=6$, and $p_{max}=24$, so the rankability given by (\ref{eq:Re}) for this graph is \[R_{e}(G)=1-\frac{1(2)}{6(24)}\approx 0.986.\]

\subsection{Spectral Rankability}
The authors of \cite{specR} introduced an alternative to $R_{e}(G)$. Referred to as spectral rankability and denoted by $R_{s}(G)$, this is a rankability measure based on the spectral variation of the Laplacian matrix of a weighted directed graph from that of a complete dominance graph.  Before discussing this, we need to recall some definitions.  

Consider a directed, weighted, graph $G$ with vertices $\{1, \ldots, n\}$, where $w_{ij}$ denotes the weight of the edge $(i,j)$.  We assume that $0 \leq w_{ij} \leq 1$ for all $i,j$.  The out-degree of vertex $i$ is then $d^+_i = \sum_{j=1}^n w_{ij}$ and the out-degree matrix $D$ is the diagonal matrix whose $i$th diagonal entry is $d^+_i$.  If $A$ is the weighted adjacency matrix of $G$, (with entries given by $w_{ij}$), the Laplacian of $G$ is then $L = D-A$.  In keeping with standard linear algebraic notation, $\sigma(M)$ denotes the spectrum (set of eigenvalues) of a square matrix $M$.

The definition of spectral rankability relies on the Hausdorff distance between two sets.  For finite sets of complex numbers $A$, $B$, we first define the quantity:
\[\sv_B(A) := \max_{x\in A} \min_{y \in B} |x-y|.\]
The Hausdorff distance between $A$, $B$ is then given by

\begin{equation}\label{eq:hd}
\hd(A,B) = \max\{\sv_B(A), \sv_A(B)\}.
\end{equation}
For two matrices $L$, $M$, we will denote the Hausdorff distance between their spectra by $\hd(L,M)$. 

Now with the Hausdorff distance defined, the spectral rankability of a directed graph $G$ with Laplacian matrix $L$ is defined as

\begin{equation}\label{eq:Rs}
R_{s}(G)=1-\frac{\hd(D, S)+\hd(L, S)}{2(n-1)},
\end{equation}
where $S$ is a diagonal matrix with $S_{ii}$ equal to $n-i$. Using this definition, the spectral rankability of graph $G$ in Figure ~\ref{Example graph for calculating edgeR} is 

\[R_{s}(G)=1-\frac{1+1}{2(4-1)}=\frac{2}{3}.\]

\subsection{Spectral and edge rankabilty of cycles: mathematical results}
\label{maths}
In the next two results, we note the behaviour of the two rankability measures on cycles.  As mentioned above, intuitively we expect that a cycle has low rankability.  It is interesting to note that the following results show that $R_{s}(G)$ and $R_{e}(G)$ behave surprisingly differently with respect to cycles, particularly as the length of the cycle grows. 

\begin{proposition}\label{prop:edgeRcycle}
Let $C$ be a cycle of length $n$.  The edge rankability of $C$, denoted $R_{e}(C)$ is given by

\begin{equation}\label{eq:edRcycle} 
R_{e}(C) = 1-\frac{2+(n-1)(n-2)}{n!(n-1)}.
\end{equation}

\end{proposition}

\textbf{Proof:} It is not hard to see that in order to transform a cycle of length $n$ to a complete dominance graph, we must:

\begin{itemize}
    \item remove at least 1 edge to obtain an acyclic graph;
    \item add at least 
    \[n-2+n-3+\cdots + 1 = \frac{(n-2)(n-1)}{2}\]
    edges.
\end{itemize}

Therefore the number of edge additions and removals required to obtain a complete dominance graph is at least $1+\frac{(n-2)(n-1)}{2}$.  Furthermore, it is not hard to see that we can obtain exactly $n$ distinct complete dominance graphs with this number of additions/removals: one corresponding to each choice of edge to remove.  It follows that in the definition of $R_{e}(G)$, $k=1+\frac{(n-2)(n-1)}{2}$, and $p=n$. This implies that

\begin{eqnarray*}
R_{e}(C) &=& 1- \frac{n(2+(n-2)(n-1))}{n!n(n-1)}\\
&=& 1- \frac{2+(n-1)(n-2)}{n!(n-1)}
\end{eqnarray*}
as claimed. 

\vspace{3mm}

\textbf{Remark:}  Note that it follows from the previous proposition that as $n$ tends to $\infty$, the edge rankability of a cycle of length $n$ tends to 1.  This is counter to the intuition that cycles should have low rankability.  In the next result, we derive a formula for the spectral rankability of a cycle of length $n$. 

\begin{proposition}\label{prop:specRcycle}
Let $C$ be a cycle of length $n$ where each edge has weight 1.  The spectral rankability of $C$, denoted $R_{s}(C)$ is given by

\begin{equation}\label{eq:spRcycle} 
R_{s}(C) = \begin{cases}
1-\frac{2n-5}{2n-2} & \mbox{$n$  even}\\
1-\frac{n-2+|n-2-e^{-i\pi/n}|}{2n-2} & \mbox{$n$ odd}.
\end{cases}
\end{equation}

\end{proposition}

\textbf{Proof:} First note that for a cycle of length $n$, the matrix of out-degrees, $D$, is simply the identity matrix as each vertex has out-degree 1.  Hence $\sigma(D) = \{1\}$ in this case and $\sigma(S)=\{0,1,\ldots, n-1\}$.  The Laplacian $L$ of a cycle is given by $I-A$ where $A$ is a cyclic matrix with $a_{i,i+1}=1$ for $1\leq i \leq n-1$, $a_{n1} = 1$, and $a_{ij} = 0$ otherwise.  It is well known that the eigenvalues of $A$ are given by the $n$th roots of unity and hence, as $L=I-A$

\[\sigma(L) = \{1-e^{2k\pi i/n}: 0\leq k \leq n-1\}.\]

It is straightforward to see that $\sv_S(D) = 0$, $\sv_D(S)=n-2$ so $\hd(S,D) = n-2$.  

Next note that $\sv_S(L) \leq 1$ as the elements of $\sigma(L)$ all lie on a circle of unit radius centred at 1.  The calculation of $\sv_L(S)$ requires a little more care.  It is clear that $\sv_L(S) = \min_{y \in \sigma(L)} |n-1-y|$.  If $n$ is even, $2 \in \sigma(L)$ and hence $\sv_{L}(S) = n-3$.  For $n$ odd, the elements of $\sigma(L)$ that are closest to $n-1$ are given by 

\[1-e^{\pi i - \frac{\pi i}{n}} = 1+e^{- \frac{\pi i}{n}},\quad 1 - e^{\pi i + \frac{\pi i}{n}}=1+e^{\frac{\pi i}{n}}\]
corresponding to $k=\frac{n-1}{2}$ and $k=\frac{n+1}{2}$ respectively.  It follows that for $n$ odd, 

\[\sv_S(L) = |n-2-e^{- \frac{\pi i}{n}}|\]
and hence (as $n \geq 2$) $\hd(S,L) = |n-2-e^{- \frac{\pi i}{n}}|$ for $n$ odd, while our earlier calculation shows that for $n$ even $\hd(S,L)=n-3$.  Putting all of this together we see that 

\[R_{s}(C) = \begin{cases}
1-\frac{2n-5}{2n-2} & \mbox{$n$  even}\\
1-\frac{n-2+|n-2-e^{-i\pi/n}|}{2n-2} & \mbox{$n$  odd}
\end{cases}\]
as claimed. 

\vspace{3mm}
\textbf{Remark:} Note that, in contrast to $R_{e}(C)$, as $n$ tends to $\infty$, the spectral rankability $R_{s}(C)$ of a cycle of length $n$ tends to 0. 

\subsection{Comparison of edge rankability and spectral rankability}
One advantage of $R_{e}(G)$ is that it is intuitive and easily understandable. A major disadvantage, however, is the computational cost of using this measure on large datasets. This is because $k$ can be calculated using an integer programme but to find an exact value for $p$ it is necessary to consider each of the $n!$ possible orderings of the vertices of the graph. This process can be completed more efficiently using a method called ``parallel exhaustive enumeration with pruning'' \cite{edgeR}, but an exact solution can still take hours to compute.

There is a reward for considering every possible ordering of the vertices, however. By examining all complete dominance graphs that are reachable by making just $k$ edge changes, it is possible to gain valuable insights into the rankability of a graph. For example, in Section ~\ref{EdgeR} we saw that there are two equally appropriate orderings of the vertices of the graph in Figure 1. A ranking algorithm might misleadingly suggest only one of these two possibilities. Similarly, $R_{s}(G)$ does not consider these two possibilities either, and unlike $R_{e}(G)$, is only capable of reflecting the rankability of a graph in a real number. In this way, $R_{e}(G)$ provides more functionality than $R_{s}(G)$, and in their paper, the authors of \cite{edgeR} describe a variety of analyses that can be performed with $R_{e}(G)$ when investigating the rankability of a dataset.

The biggest advantage of $R_{s}(G)$ over $R_{e}(G)$ is that it is much more widely applicable. This is made possible by the existence of efficient algorithms for the computation of matrix eigenvalues such as the QR algorithm \cite{qr}. A disadvantage, however, is that the use of eigenvalues rather than edges makes this method less intuitive. A second disadvantage is that owing to the comparison with a complete dominance graph, this measure struggles to return meaningful results where few of the pairwise comparisons have been considered, such as in the case of sparse data.

\subsection{Weighted Edges}
Our description of the rankability problem thus far has mainly considered the unweighted case, i.e. data consisting of binary relationships between items; item $i$ beat item $j$, represented by a 1 in position $(i,j)$ of the adjacency matrix, or item $i$ did not beat item $j$, represented by a 0 in position $(i,j)$ of the adjacency matrix. Often, relationships between items within a dataset are more complex than this, and the amount by which one item has beaten another may be of interest. To capture this extra information regarding the pairwise comparisons between a set of items it is necessary to use weighted edges. For example, instead of ranking film $i$ more popular than film $j$ because the majority of people surveyed gave it a higher rating on Netflix, a weighted adjacency matrix could represent this difference in popularity by assigning to the entry in position $(i,j)$ of the adjacency matrix the proportion of people that having watched both film $i$ and film $j$, decided they liked film $i$ more.

Recently in \cite{weighted_rankability}, Anderson et al. investigated this weighted case of rankability. In particular, they extend the notion of \emph{distance from perfection} to weighted data. They also show how indirect comparisons between items within a dataset can be used to mitigate the challenges posed by assigning a rankability to sparse data. This approach involves adding weighted edges from an item $i$ to an item $j$ if item $i$ has beaten a third item $k$ which has demonstrated dominance over item $j$ in a direct comparison.

In this paper, however, we will focus on the unweighted case of rankability, and will tackle the challenge posed by sparse data with a method which penalises sparse datasets with a lower rankability. Our approach, instead of attempting to remove sparsity from a dataset, accepts sparsity as an innate feature of the dataset, and adjusts the rankability of such datasets accordingly. The absence of a mechanism for incorporating weighted edges directly into the original graph or adjacency matrix when using our supervised learning approach is perhaps only a slight disadvantage of our method. Indeed, many ranking algorithms actually choose to avoid weighted edges in favour of unweighted edges \cite{colley}.

\section{Rankability Measure Evaluation}
Our results from Section~\ref{maths} demonstrate the appeal of being able to investigate trends in rankability for families of graphs with known properties. This enabled us to compare the rankability from $R_{e}(G)$ and $R_{s}(G)$ with our own intuition which identifies cyclic graphs as being inherently unrankable. It would also be desirable to extend this analysis to other families, but it may not always be possible to easily derive an expression for the rankability of other types of graphs.

With the above observations in mind, we now propose a method for studying the rankability of data in a situation where we are able to generate synthetic data/graphs with an estimate of their rankability, which we call the \emph{target rankability}. The method for generating the synthetic data and the target rankability could be chosen in a number of ways, and could be application or data dependent. As our main example, we will generate the synthetic graphs by making random modifications to the complete dominance graph. The larger the probability of modification used in generating the synthetic data, the lower the assigned target rankability.

In this section, to demonstrate that this is a feasible approach, we show that this approach can generate graphs with target rankabilities that are usefully correlated with the rankability metrics introduced in Section~\ref{sec:previous}. We can think of this as either validating the previous rankability metrics or as providing evidence for the usefulness of the target rankability metric chosen.

\subsection{Target Rankability}
The existing rankability measures discussed in the previous sections  measure rankability in different ways, and both will usually return a different result, even when applied to the same dataset. In the case of small graphs, it is possible to scan by eye the edges that exist within a graph and form an intuitive idea about its rankability. In this way it is possible to verify that the rankability measures are returning a sensible result. But how can we check if a rankability measure is performing well when there are 20, 50, 100, even 1000 or more vertices? We now seek to develop the framework necessary to answer this question.

Our approach will be to define a graph model to generate graphs. To each graph generated by this model we will assign a target rankability determined by the parameter(s) of the model that created it. By applying different rankability measures to these graphs it should be possible to determine the extent to which these measures are correlated with the target rankability that we have selected. We emphasise that the following definition of target rankability is an example. There are potentially many ways of generating graphs with a target rankability selected to match a particular application, of which this method is just one.

We will now describe our example model for generating graphs. Consider the adjacency matrix $A$ of a 4 vertex complete dominance graph shown below. 

\[A=
\begin{blockarray}{ccccc}
&1&2&3&4\\
\begin{block}{c(cccc)}
  1&0 & 1 & 1 & 1 &\\
  2&0 & 0 & 1 & 1 &\\
  3&0 & 0 & 0& 1 &\\
  4&0 & 0 & 0 & 0 &\\
\end{block}
\end{blockarray}
 \] 
 
We will say that ``switching'' the status of an edge at position $A_{ij}$ has the effect of turning a 1 into a 0 and vice versa. Our probabilistic graph generation model is as follows:

\begin{enumerate}
\item Start with a complete dominance graph and a parameter $0\leq p\leq 0.5$.
\item Switch the status of every off-diagonal entry in the adjacency matrix with probability $p$.
\end{enumerate}

The reason for choosing a value for $p$ between 0 and 0.5 is that as the value of $p$ approaches 1, the effect of this is to produce a new complete dominance graph with the direction of every edge reversed. If the value of $p$ is 0 then the adjacency matrix will be completely unchanged. If the value of $p$ is 0.5 then the effect of this is to produce a completely random graph where every edge exists (or doesn't exist) with probability $p=0.5$. Thus, for $0 \le p \le 0.5$, larger $p$ values will tend to represent graphs that are less rankable. To convert the parameter $p$ into a number that can be interpreted as a measure of the graph's rankability we choose to define the target rankability of a graph produced with this model as

\begin{equation}\label{eq:target1}
t=1-2p.
\end{equation}

Initially it might appear strange that we have made the decision to define the target rankability of a graph in terms of a complete dominance graph. This decision might make it difficult to obtain meaningful results when measuring the rankability of sparse data. Later we will make a slight modification to (\ref{eq:target1}) to account for this, but for now we will only work with complete data.

\subsection{Diagnostic Plots}
\begin{figure}[ht]
    \includegraphics[width=\linewidth]{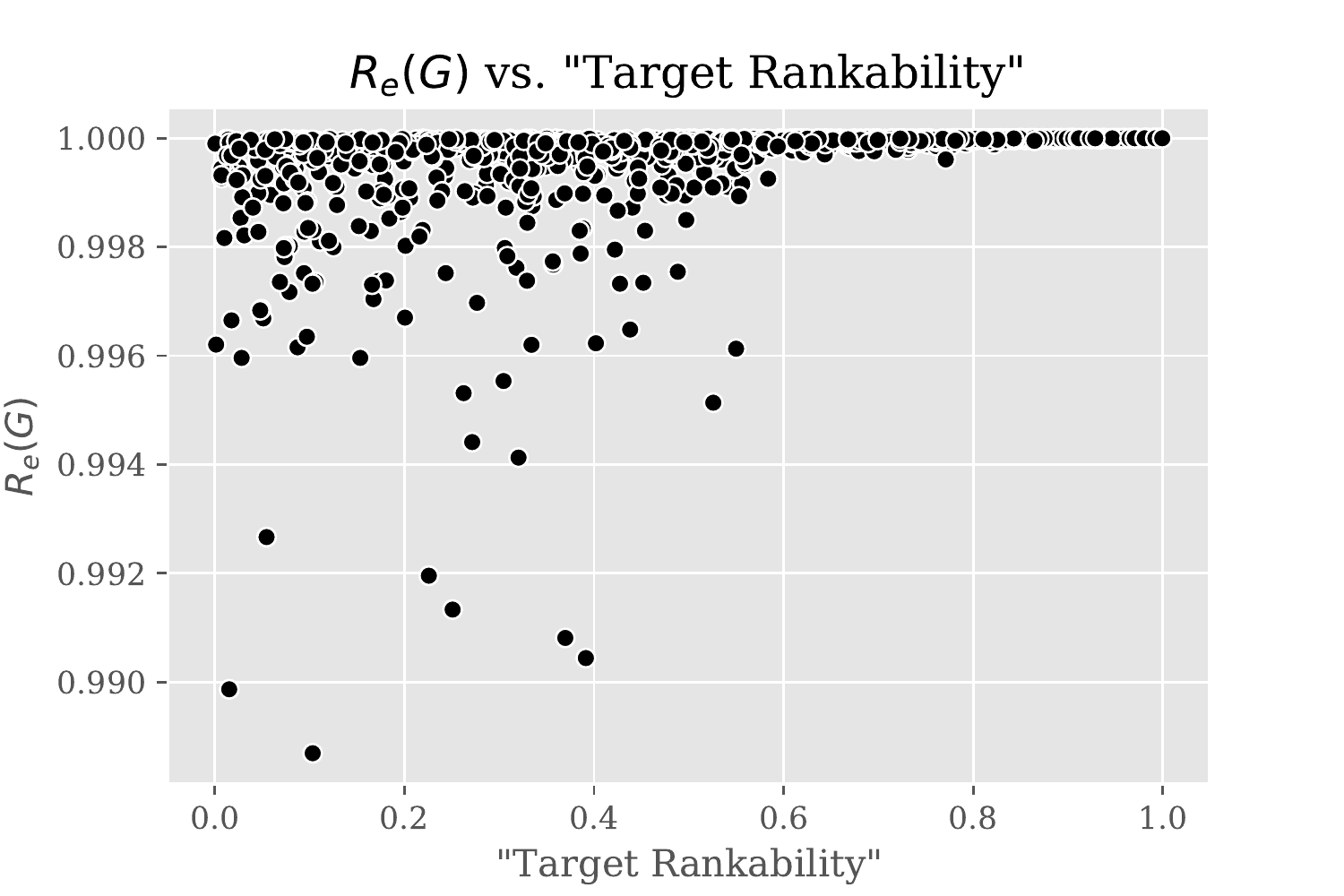}
    \caption{$R_{e}(G)$ vs. target rankability, $n=8$ simulation, $\rho=0.747$}
    \label{fig:Revst}
\end{figure}

\begin{figure}
    \includegraphics[width=\linewidth]{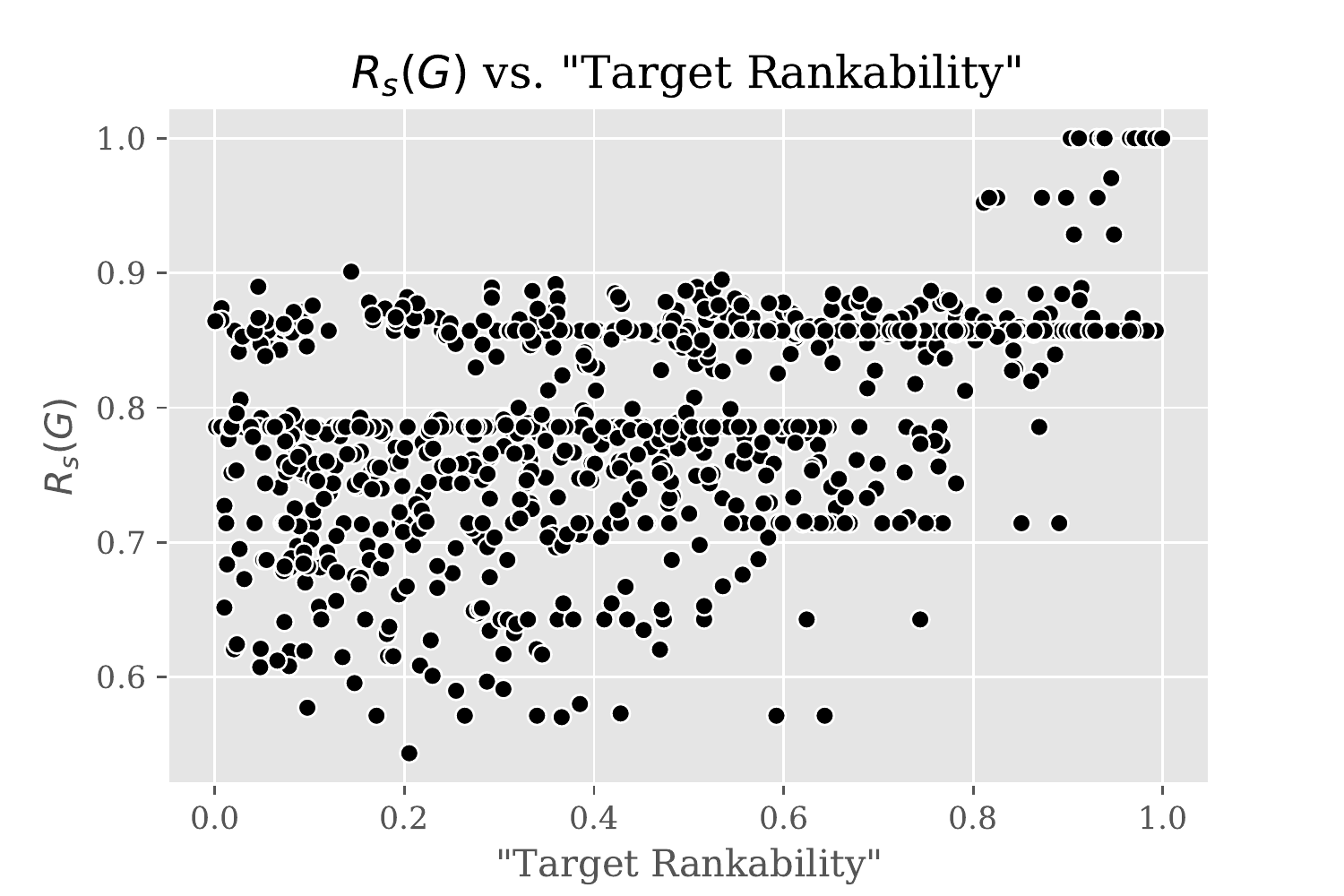}
    \caption{$R_{s}(G)$ vs. target rankability, $n=8$ simulation, $\rho=0.584$}
    \label{fig:Rsvst}
\end{figure}

A graph's target rankability should provide us with a good indication of the rankability of the graph. By applying the rankability measures $R_{e}(G)$ and $R_{s}(G)$ to graphs produced with the graph generation process described above, it should be possible to determine the extent to which the rankability measures are associated with the target rankability of such graphs. In the following examples we have completed 1000 simulations, each with a randomly selected value for the parameter $p$ in (\ref{eq:target1}) between 0 and 0.5. Figure~\ref{fig:Revst} shows the association between the rankability measure $R_{e}(G)$ and the target rankability of  graphs with $n=8$ vertices. The results from $R_{e}(G)$ are consistently very close to 1 because of the large denominator in (\ref{eq:Re}), even when the target rankability is 0. The results from $R_{s}(G)$ in Figure~\ref{fig:Rsvst} appear slightly better, but are quite quantised along the $R_{s}(G)$ axis. 


For comparison, we have also included $\rho$, the Spearman rank correlation coefficient, correct to three decimal places in the description for each graph, and in Table~\ref{tab:spearman}. The coefficients are all relatively close to 1, indicating that the target rankability is ordering the rankability of the graphs in a similar way to previous metrics. Interestingly, the Spearman rank correlation coefficient for the relationship between $R_{e}(G)$ and the target rankability is greater than that of the relationship between $R_{s}(G)$ and the target rankability. This relationship is only detectable by analysing $R_{e}(G)$ differences of less than the order of $10^{-3}$. Such small differences are not easily interpretable when examining the rankability of a single graph. Instead, these results suggest that the strength of $R_{e}(G)$ lies not in the rankability that it assigns to a graph, but in its ability to analyse each of the $n!$ orderings of its vertices as we discussed in Section~\ref{sec:previous}.

\begin{table*}
\begin{center}
\caption{Spearman Rank Correlations, Target Rankability}
\label{tab:spearman}
\begin{tabular}{lcccccl}\toprule
& \multicolumn{3}{c}{Complete Data} & \multicolumn{3}{c}{Sparse Data}
\\\cmidrule(lr){2-4}\cmidrule(lr){5-7}
           & $n=8$  & $n=20$ & $n=50$ & $n=8$  & $n=20$ & $n=50$\\\midrule
$R_{e}(G)$    & 0.747 & NA & NA & 0.540 & NA & NA \\
$R_{s}(G)$ & 0.584 & 0.888 & 0.969 & 0.462 & 0.657 & 0.720\\
$R_{f}(G)$   & 0.868 & 0.968 & 0.994 & 0.676 & 0.912 & 0.976 \\\bottomrule
\end{tabular}
\end{center}
\end{table*}

Due to the computational cost of running $R_{e}(G)$ on larger datasets it is not possible to make any comparisons on larger graphs. Table~\ref{tab:spearman} includes the results from applying $R_{s}(G)$ to graphs with $n=20$ and $n=50$ vertices respectively. With larger graphs, the results from $R_{s}(G)$ appear to be much better than before and there is a  positive correlation with the target rankability, that appears to be increasing with $n$.

$R_{e}(G)$ and $R_{s}(G)$ both rely on a comparison with a complete dominance graph. The authors of \cite{specR} identify this as a potential limitation of $R_{s}(G)$, because this comparison with a complete dominance graph can make it difficult to obtain meaningful results when measuring the rankability of sparse data. To explore this we have prepared graphs with a known target rankability but have then removed some of their edges. Specifically, we have generated 1000 graphs with $n$ vertices, all generated with a random parameter $p$ in (\ref{eq:target1}). After creating each of these graphs we have then generated a second random parameter $c$, and have removed every edge in the graph with probability $1-c$. To take this effect into account we have defined the target rankability of such graphs as $t=c(1-2p)$. The intuition behind this decision is that as fewer of the pairwise comparisons remain, the less rankable the dataset is likely to be. In the presence of sparse data, it is clear that the relationship between $R_{e}(G)$ and $R_{s}(G)$ with the target rankability can become quite weak. This can be seen by examining the correlations with sparse data contained in Table~\ref{tab:spearman}. Again, $R_{s}(G)$ appears to perform slightly better on the larger graphs, and, as expected, there is a weaker correlation with the target rankability when compared to the results for graphs which contain no missing data.

\section{Regression Using Target Rankability}
In this section we will describe how, using our model to generate graphs and our target rankability, we can extend these to provide a rankability measure based on the properties of a graph or its associated adjacency matrix in a supervised learning approach to rankability.

\subsection{Training a Rankability Model}
The simulations used in the previous section enable us to reframe the task of rankability as a supervised learning problem. The idea is to use the target rankability as the response variable and predict it via properties of a graph or its adjacency matrix to measure its rankability. The system can then assign a rankability to any graph, without knowing the model that generated it.

To predict the target rankability, we suggest using graph properties that are easy to calculate, while also having relevance to ranking. The properties we will use as predictor variables are:

\begin{enumerate}
\item The number of triangles in a graph. We will only count unique triangles, and will say that the triangle $i\rightarrow j\rightarrow k\rightarrow i$ is the same as $j\rightarrow k\rightarrow i\rightarrow j$. 
\item The number of 2-cycles, which we will call contradictions, because they represent cases where $i$ has beaten $j$ but $j$ has also beaten $i$. 
\item The standard deviation of the out-degrees of the vertices of the graph.
\item The directed algebraic connectivity of a graph as defined by Wu in \cite{wu}.
\item The number of times two vertices in a graph are not directly connected. We will refer to such situations as draws because there is no evidence to suggest either vertex $i$ or $j$ is better than the other.
\end{enumerate}

We have chosen to use these properties for a number of reasons. Firstly, in the case of the number of triangles in a graph, our intuition is that graphs with fewer triangles are likely to be more rankable because cycles are inherently difficult to create an ordering from. Similarly, we also expect a graph with fewer 2-cycles will also be more rankable.
The out-degree counts the number of out-going edges from a given vertex. We expect that unrankable graphs will have no clear ``winners'' or ``losers'', meaning many vertices will have a similar number of ``wins''. Consequently, we expect that unrankable graphs will have a smaller standard deviation of the out-degrees relative to graphs which are more rankable. 
With regards to the directed algebraic connectivity of a graph, we expect graphs with higher connectivities will be less rankable for the same reasons as having a higher number of 2-cycles or 3-cycles. Directed algebraic connectivity was also suggested as a helpful property for quantifying rankability by the authors of \cite{specR}.  
Finally, graphs with a large number of draws are missing large amounts of data. Without a direct comparison between two vertices it is very difficult to rank one higher than the other. For this reason, we expect graphs with a large number of draws to be less rankable.

We also note that none of these properties are prohibitively computationally expensive to calculate. This means there is no need for us to restrict our approach to graphs of a relatively small size, and makes our approach widely applicable to many problems.

Our first step will be to train a random forest regression model on synthetic data generated by the target-rankability graph generation process. We use the python package scikit-learn with no modifications to the default tuning parameters \cite{scikit-learn}. There are two main advantages of using a random forest here. Firstly, a random forest will not make predictions outside of the range of the original training values. This can be a drawback of using random forests, but in this context it is advantageous because we only want to make predictions between 0 and 1. Secondly, in the case of sparse data, the association between the target rankability and the predictor variables can become strongly non linear, making a linear model less appropriate. Other modelling approaches are certainly possible however.

We will denote the rankability returned by this model as $R_{f}(G)$. To evaluate the efficacy of this approach we have trained a random forest on 1000 graphs with $n$ vertices, all generated with a random parameter $p$ in (\ref{eq:target1}). We have then tested the random forest on 1000 different $n$ vertex graphs generated by the same process. Finally, we have calculated the Spearman rank correlation coefficient to assess the performance of the model, as we did with the other rankability measures in the previous section, and have added our results to Table~\ref{tab:spearman}.
To assess the performance of $R_{f}(G)$ in the case of sparse data, we have trained and tested the random forest regression model on graphs with target rankability equal to $t=c(1-2p)$ using the procedure described previously.

In both cases, the correlation between the predicted rankability $R_{f}(G)$, and the target rankability of different sized graphs in Table~\ref{tab:spearman} is strong.
This is promising, but our tests are thus far limited to the case where we have trained and tested our random forest on graphs that have both been generated with the same graph generation model. To improve our confidence that a regression approach to rankability is applicable to other types of graphs we use a second, independent, graph generation model for testing our trained random forest. We will use the graph generation technique for simulating a competition between a set of $n$ items described in \cite{graphModel}. We will also need a means to determine the relative rankability of the graphs generated using this technique. This is the aspect that we will address first.

\subsection{Validation with Out of Sample Synthetic Data}
\begin{table*}
\begin{center}
\caption{Spearman Rank Correlations, Relative Rankability}
\label{tab:spearman2}
\begin{tabular}{lcccccl}\toprule
& \multicolumn{3}{c}{Complete Data} & \multicolumn{3}{c}{Sparse Data}
\\\cmidrule(lr){2-4}\cmidrule(lr){5-7}
           & $n=8$  & $n=20$ & $n=50$ & $n=8$  & $n=20$ & $n=50$\\\midrule
$R_{e}(G)$    & 0.667 & NA & NA & 0.564 & NA & NA \\
$R_{s}(G)$ & 0.726 & 0.931 & 0.969 & 0.477 & 0.614 & 0.680\\
$R_{f}(G)$   & 0.882 & 0.979 & 0.995 & 0.764 & 0.911 & 0.967 \\\bottomrule
\end{tabular}
\end{center}
\end{table*}

Consider a set of $n$ items each with an ability level $a_{i}$. Using the Elo rating system \cite{elo}, which was originally invented for rating chess players, the expected outcome for item $i$ against item $j$ is given by \[E_{ij}=\frac{1}{1+10^\frac{a_{j}-a_{i}}{400}},\] and $E_{ji}$ is given by $1-E_{ij}$. To enable the exact calculation of the probability of the pairwise comparison between item $i$ and item $j$ resulting in a draw we will make the following modification: Define the probability of a draw as $P_{d}=E_{ij}E_{ji}$. This will give a higher probability of a draw as the two items get closer in ability levels. Now define the probability that item $i$ will beat item $j$ as $P_{ij}=E_{ij}(1-P_{d})$. Finally define the ``relative rankability'' of a set of $n$ items as 

\begin{equation}\label{eq:relrank}
r=\frac{2\sum_{i<j}|P_{ij}-P_{ji}|}{n(n-1)},
\end{equation}
which is the average of the absolute values of $P_{ij}-P_{ji}$ over all cases where $i$ is not equal to $j$.

Now to generate a graph with a known ``relative rankability'' we use the following procedure:

\begin{enumerate}
\item Begin with a completely disconnected graph with $n$ vertices.
\item Assign to each vertex an ability level $a_{i}$. 
\item For each pair $(i, j)$, calculate the probabilities $P_{d}$, $P_{ij}$, and $P_{ji}$.
\item Calculate the relative rankability of the graph using (\ref{eq:relrank}).
\item For each pair $(i, j)$, generate a uniformly distributed random number $0\leq \tau \leq1$. If $\tau<P_{ij}$ then add an edge from vertex $i$ to vertex $j$. If $\tau>P_{ij}+P_{d}$ then add an edge from vertex $j$ to vertex $i$.
\item Optionally repeat step 5 above to allow for the possibility of 2-cycles in the generated graph.
\end{enumerate}

As we mentioned earlier, there are potentially many ways of creating a graph with a ``known'' rankability, and this method which we refer to as ``relative rankability'' is a member of this much larger family of methods.

In order to evaluate the efficacy of $R_{f}(G)$ on out of sample synthetic data we trained a random forest on 1000 graphs with $n$ vertices, all generated with a random parameter $p$ related to target rankability as in (\ref{eq:target1}). We have then tested the random forest on 1000 different $n$ vertex graphs generated with a known ``relative rankability'' using the procedure above (including step 6 to generate possible ties). Finally, we calculate the Spearman rank correlation coefficient to assess the performance of the model. These results can be found in Table~\ref{tab:spearman2}. The predictions made by the random forest still show a clear increasing relationship with the other rankability metrics, that appears to strengthen as $n$ increases. This suggests that the random forest regression model is performing well, even on out-of-sample synthetic data. To obtain the results for sparse data in Table~\ref{tab:spearman2}, we have trained the random forest regression model on graphs generated using the sparse target-rankability procedure and $t=c(1-2p)$. We have then tested the trained model on graphs generated using the ``relative rankability'' procedure which have had edges removed with probability $c$, and
\begin{equation}
r=\frac{2c\sum_{i<j}|P_{ij}-P_{ji}|}{n(n-1)}.
\end{equation}

\section{Case Study: Non-Synthetic Data}

\begin{figure}
\begin{center}
\includegraphics[width=\linewidth]{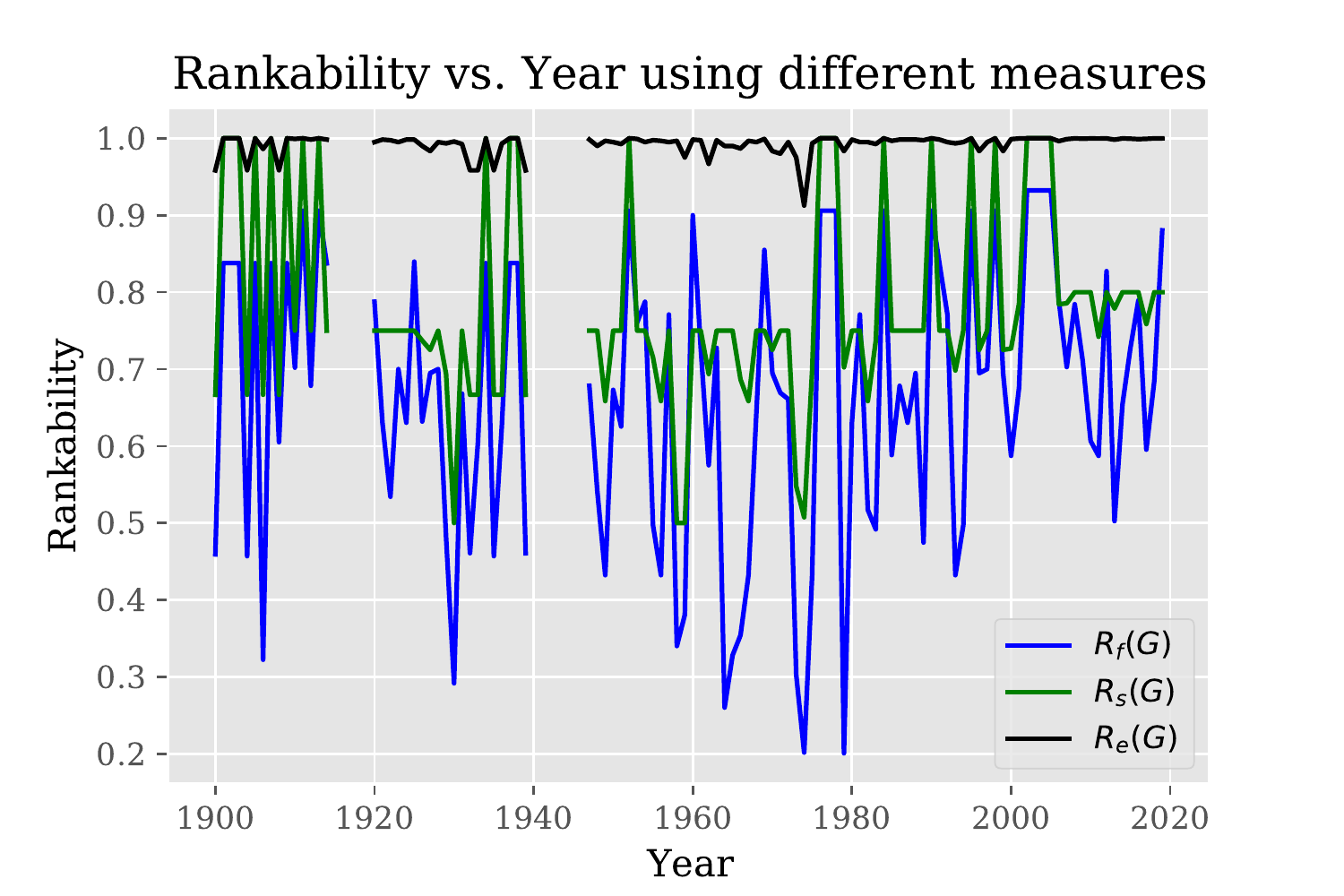}
\captionof{figure}{Rankability of rugby competitions}
\label{Rankability of rugby competitions}
\end{center}
\end{figure}

In this section we will apply the three rankability measures, $R_{e}(G)$, $R_{s}(G)$, and $R_{f}(G)$ to real-life sports data from rugby and football competitions. Our first dataset contains information on all rugby games played in the Home Nations, Five Nations, and Six Nations Championships from 1900 to 2019 inclusive \cite{rugbyData}. The Home Nations Championship is a four-team event, and the Five and Six Nations Championships are five- and six-team events. In these competitions, each team plays every other team once. To construct a graph for every year we have started with a completely disconnected graph and added an edge from vertex $i$ to vertex $j$ if and only if team $i$ has beaten team $j$ in the game involving these two teams. The rankabilities of each year's graph using the different rankability measures are shown in Figure~\ref{Rankability of rugby competitions}. When the competition involved four teams we have used a random forest that was trained on graphs with four vertices, and when the competition involved five or six teams we have used a random forest that was trained on graphs with five or six vertices respectively.

\begin{figure*}
\begin{center}
\includegraphics[width=\linewidth]{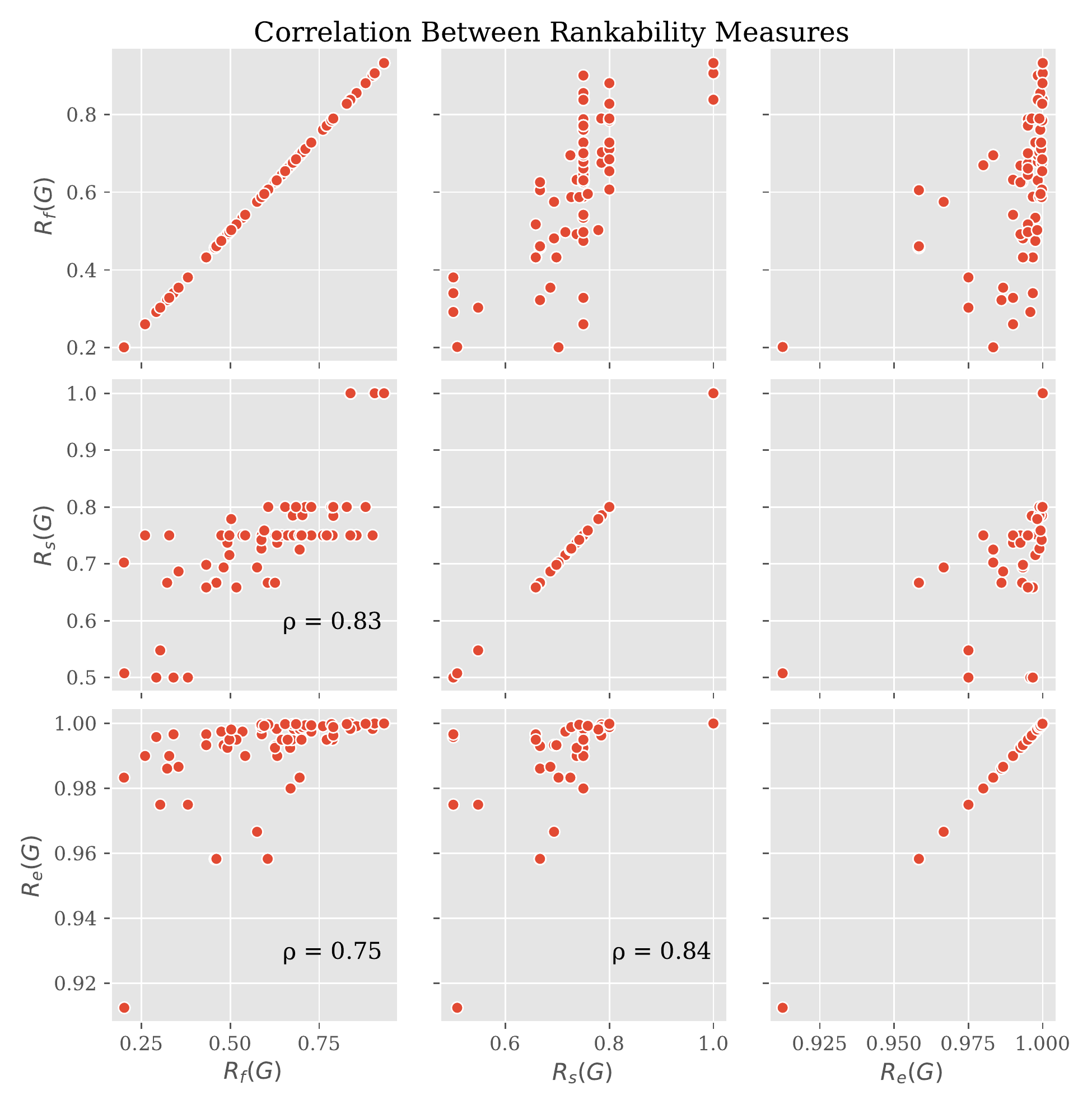}
\captionof{figure}{Rankability of Rugby Competitions Correlations}
\label{Rugby Correlations}
\end{center}
\end{figure*}

Unlike the other rankability measures, the random forest fails to return a rankability of 1 for the years where the graph was a complete dominance graph. The reason for this is that the random forest is trained on graphs which are determined by a parameter $p$ in (\ref{eq:target1}), and a low value for $p$ close to 0 will still occasionally produce a complete dominance graph (especially in the setting with small $n$ as is the case here). All rankability measures identify 1974 as a year with a very low rankability. This year every team lost at least one game, three games were drawn, and three of the five teams finished on the same number of points at the end of the competition. As can be seen from Figure~\ref{Rankability of rugby competitions}, all three rankability measures approximately follow the same trends and respond to the same spikes or drops in the rankability of the matches. This agreement is reflected in a strong correlation between the results. Figure~\ref{Rugby Correlations} shows a scatterplot matrix of these results and also gives the Spearman rank correlation coefficients.

\begin{figure}
\begin{center}
\includegraphics[width=\linewidth]{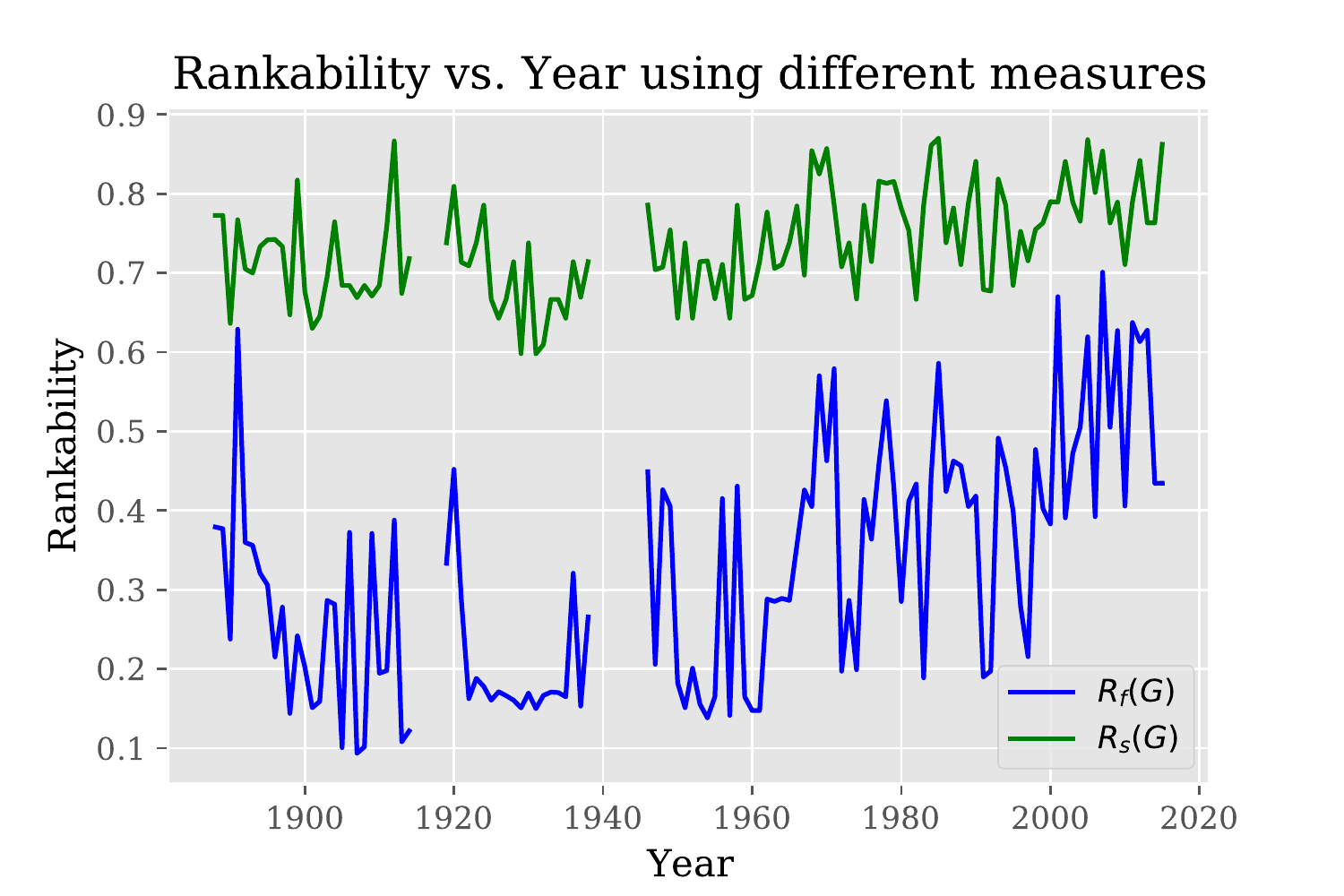}
\captionof{figure}{Rankability of football competitions}\label{Rankability of football competitions}
\end{center}
\end{figure}

Next, we turn our attention to a dataset containing every game played in the top-tier of English football from 1888 to 2015 inclusive \cite{footballData}. In this competition, every team plays every other team twice. We have constructed the graph for each season by starting with a completely disconnected graph and adding an edge from vertex $i$ to vertex $j$ if and only if team $i$ beat team $j$ at least once. For each season, we have used a random forest  that was trained on graphs with the same number of vertices as the number of teams in that season. Figure~\ref{Rankability of football competitions} shows the rankability of each of the seasons in this dataset. $R_{e}(G)$ was not applied to this dataset because the number of teams/vertices is too large, and this would be too computationally expensive.

\begin{figure}
\begin{center}
\includegraphics[width=\linewidth]{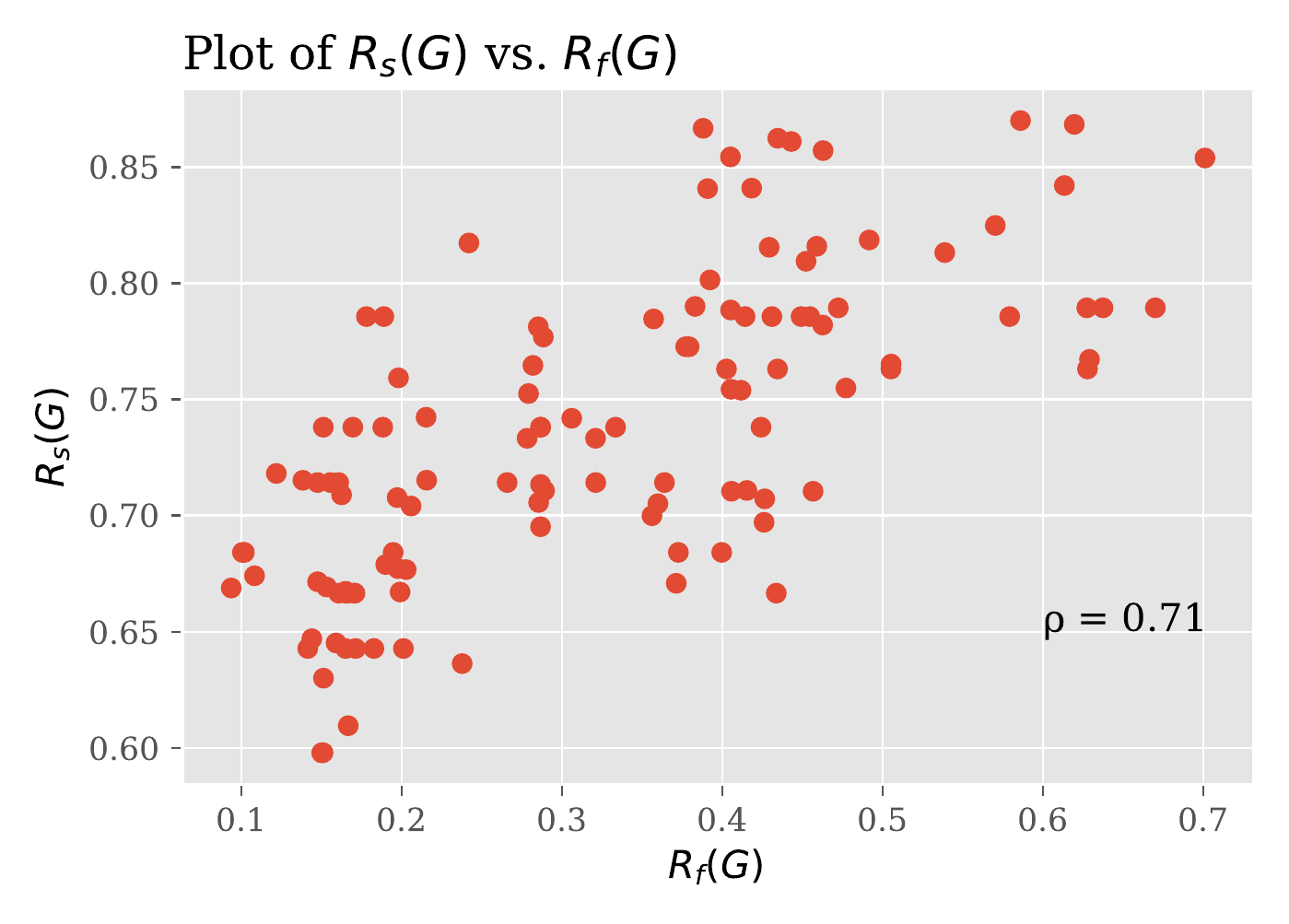}
\captionof{figure}{Rankability of Football Competitions Correlations}
\label{Football Correlations}
\end{center}
\end{figure}

Both rankability measures assign a relatively high rankability to the season beginning in 2007. This season ended with just two teams finishing on the same number of points, and this is a relatively rare occurrence in the premier league. This season also contains a relatively low number of draws, contradictions and triangles when compared to the other seasons involving 20 teams. At the opposite end of the spectrum, both rankability measures assign a relatively low rankability to the season beginning in 1957. This season contains a relatively high number of draws, contradictions and triangles when compared to the other seasons involving 22 teams. Again we can see that $R_{f}(G)$ and $R_{s}(G)$ jointly track many of the spikes and drops in the rankability of these seasons. This high level of agreement is once again reflected by the strong correlation between the results from these rankability measures, shown in Figure~\ref{Football Correlations}.

If there was no correlation between the rankability measures in this section then it might have suggested that they were performing poorly because they could not agree on what the most and least rankable seasons were. The fact that there is a moderately strong correlation between them, however, suggests that they could be performing well because they are in approximate agreement regarding the rankability of the seasons.

\section{Discussion}

In this paper we have developed a framework for rankability based on the generation of graphs with a chosen/candidate \emph{target rankability}. The underlying idea involves making perturbations to complete dominance graphs. With this approach we have been able to achieve two things. Firstly, we have shown that our approach enables evaluation of the performance of different rankability measures. This is important because without a method for the evaluation of a rankability measure’s performance, it is difficult to determine how effective it is on larger graphs. 

This approach is very useful because it allows us to visualise how the rankability measures respond as the graphs transition from being completely random to a complete dominance graph. It is also very intuitive because it allows us to easily interpret the rankability of these graphs in terms of the parameters used in the probabilistic graph model which created them.

A disadvantage of this method for evaluating a rankability measure’s performance, however, is that it can only be used to evaluate a rankability measure on graphs that have been generated by the model (in our example, making perturbations to a complete dominance graph). Therefore, this method might not provide a good indication of how rankability measures will perform when applied to other types of graphs. However, it could clearly be extended to include other families of graphs (say, perturbed cycles).

A second key contribution of this paper describes how to use a supervised learning approach for the measurement of rankability. We have demonstrated this approach by training a random forest on synthetic data generated with a chosen target rankability. We have also tested random forest regression models on graphs generated with a different graph generation model described in \cite{graphModel}. Results from these experiments suggest that the random forest is capable of performing well on out of sample synthetic data. Further experiments show that a modification to the target rankability graph generation process can make the random forest effective in the case of sparse data also. Our experiments with real life sports data from rugby and football have shown that in many cases there is a strong positive correlation between the results from the different rankability measures. This agreement between the rankability measures also provides support to the supervised learning approach.

One of the main strengths of the supervised learning approach is that it is flexible. There are many combinations of regression models and graph properties that could be used with this approach. There are also potentially many ways of generating graphs with a ``known'' rankability. This flexibility makes it possible for researchers to tailor this approach to the dataset whose rankability they want to investigate.

We suggest that the steps of this procedure would include:
\begin{enumerate}
\item Defining an alternative graph generation model and choice of target rankability that they believe will be well suited to their data.
\item Evaluating how the existing rankability measures perform with this choice of target rankability.
\item Training their choice of regression model.
\item Applying the chosen rankability measures and the trained regression model to their data.
\end{enumerate}

One disadvantage of the supervised learning approach is that predictor variables, such as the number of triangles, will remain unchanged in graphs with weighted edges. In the case of weighted graphs, a choice of predictor variables that provide a useful representation of the weights would be necessary. Furthermore, as is the case when evaluating a rankability measure's performance, this method is based on graphs that have been generated by making perturbations to complete dominance graphs. Therefore, the training data upon which the supervised machine learning model learns will contain few graphs which are unlikely to be produced by this process such as cyclic graphs. This means that the regression models may perform poorly on rarely seen graphs such as these.

Areas for future research include the development of graph generation models with known target rankability that can account for a wider variety of graphs such as perturbed cyclic graphs. Future research could also include the identification of different graph or adjacency matrix properties that can aid the supervised learning approach. Additionally, a further generalisation of the regression approach to rankability by including the number of nodes in a graph as a predictor variable would also be of interest. This would allow the same trained model to be applied to multiple graphs of different sizes, and increase the flexibility of this approach. 

\section{Conclusion}

Rankability measures are an important tool for assessing the appropriateness of imposing a strict linear ordering on any set of items. Ideally, a rankability measure will be widely applicable, and will return sensible results which agree with our intuition based on the structure of a given graph. By working in the reverse direction, and starting with graphs of known ``target rankability", we have developed a framework capable of addressing these two requirements. This method, which is very flexible, has been demonstrated by making perturbations to complete dominance graphs, and enables the evaluation of the performance of any rankability measure. It also facilitates the use of a supervised learning approach which can ``predict" the rankability of a graph. Our results show that this method performs well both on out of sample synthetic data, and in the case of sparse data also. Results from real life sports data support our method further by highlighting a strong correlation between our approach and the existing rankability measures described in \cite{edgeR} and \cite{specR}.
 
\section*{Funding}
This work has emanated from research conducted with the financial support of Science Foundation Ireland under grant number 18/CRT/6049.

\bibliographystyle{IEEEtran}
\bibliography{main.bib}

\end{document}